\documentclass[12pt]{article}
\usepackage{latexsym}
\usepackage{amssymb}
\usepackage{mathrsfs}
\usepackage[OT2,OT1]{fontenc}
\usepackage{graphicx}
\DeclareGraphicsExtensions{.eps, .jpg}
\def\cyr{%
\renewcommand\rmdefault{wncyr}%
\renewcommand\sfdefault{wncyss}%
\renewcommand\encodingdefault{OT2}%
\normalfontx
\selectfont}
\DeclareTextFontCommand{\textcyr}{\cyr}

\newcommand{\CC}{\mathbb{C}}

\newcommand{\CP}{\mathbb C\mathbb P}
\newcommand{\RP}{\mathbb R\mathbb P}

\newcommand{\RR}{\mathbb R}

\newcommand{\ZZ}{\mathbb{Z}}

\newcommand\wee{{\mathbf w}}
\newcommand\bfv{{\mathbf v}}

\newtheorem{main}{Theorem}

\newtheorem{thm}{Theorem}

\newtheorem{lem}[thm]{Lemma}

\newenvironment{proof}{\medskip \noindent {\bf Proof.}}{\hfill $\blacksquare$\\}
\begin{document}
\title{The Einstein-Weyl Equations, \\ Scattering Maps,  and Holomorphic Disks}
\author{Claude LeBrun\thanks{Supported 
in part by  NSF grant DMS-0604735.}     ~and  
L.J. Mason\thanks{Supported in part by the EU through FP6
Marie Curie RTN {\em ENIGMA} (MRTN-CT-2004-5652).}}
\date{}
\maketitle 
%


Recent years have witnessed a resurgence  of 
twistor methods in  both mathematics and physics. 
On the physics side, this development  has primarily  
been driven by Witten's discovery  \cite{witts} 
that strings in twistor space
can be used to calculate Yang-Mills scattering amplitudes, resulting in 
concrete  experimental predictions that can be tested 
using existing  particle accelerators. Here, one of the main strands of thought
  \cite{berkwit} involves  the introduction   of {\em open strings} -----
 Riemann surfaces with non-empty boundary, where the 
 boundary is constrained to lie on a  specified submanifold,  called a $D$-brane.

On the mathematical side, there has been a parallel development. 
Penrose-type    twistor correspondences  study a component of 
moduli 
spaces of compact holomorphic curves $C$ in a complex manifold $Z$; 
when $Z$ is equipped with an anti-holomophic involution $\sigma  : Z\to Z$,
the parameter space for  those curves $C\subset Z$ with $\sigma (C) = C$  
  then often turns out to 
carry a natural geometrical structure  which is  the general solution of
an interesting differential-geometric problem.  By contrast, the new paradigm is 
instead to  study the moduli spaces of holomorphic curves-with-boundary  $C$,
where $\partial C\neq \varnothing$  is constrained to lie in a
 totally real submanifold $P\subset Z$. 
It  turns out that this framework is well adapted to the
study of many  natural problems in global differential geometry
where solutions typically are of  very low regularity. 
For example,  Zoll metrics (and, more generally, Zoll  projective structures) on surfaces 
turn out to arise \cite{lmzoll} from holomorphic disks in $\CP_2$ with boundaries
in a totally real $\RP^2\hookrightarrow \CP_2$. Analogous results \cite{lmfrei} 
describe 
split-signature self-dual conformal structures on $4$-manifolds 
in terms of holomorphic disks in $\CP_3$ with boundary on a 
totally real $\RP^3\hookrightarrow \CP_3$.  For  closely-related results on 
the Yang-Mills equations in split signature, see \cite{splitym}.

This article will show  that these techniques also lead to  definitive  results concerning 
the Einstein-Weyl equations for a $3$-dimensional Lorentzian space-time, 
thereby substantiating  a claim made in  \cite{lebrsb}.
Recall that the
Einstein-Weyl equations, for   a conformal 
class of metrics $[g]$ and a compatible  torson-free connection 
$\triangledown$,  are precisely the requirement  that the trace-free
symmetric part of the Ricci  tensor of $\triangledown$ should vanish; 
that is, the torsion-free connection $\triangledown$ is required to 
satisfy 
$\triangledown g = \alpha \otimes g$
for some $1$-form $\alpha$, and its
 curvature tensor ${\mathcal R}$  is required to satisfy 
 \begin{equation}
 \label{hew} 
{{\mathcal R}^k}_{ikj} + {{\mathcal R}^k}_{jki} = f g_{ij}
\end{equation}
for some function $f$. 
In $3$-dimensions, these  equations 
were first investigated by 
Elie Cartan \cite{cartanweyl},
who   discovered (in essence)  that they are 
totally integrable. Much later, 
Hitchin  \cite{hitproj} rediscovered  the  $3$-dimensional Einstein-Weyl equations
 as a dimensional reduction of the $4$-dimensional self-duality equations, 
and  described their  Riemannian solutions in terms of 
 the so-called mini-twistor correspondence \cite{hitmini,jonestod}. 
 The twistor correspondence we will develop here is a close cousin  of 
  Hitchin's mini-twistor correspondence, but is naturally adapted  to low-regularity 
   Lorentzian solutions
with very specific asymptotic behavior. Our main result  is the 
following:

\begin{main} 
There is a natural one-to-one correspondence between 
\begin{itemize}
\item smooth,  space-time-oriented, conformally compact,
globally hyperbolic, Lorentzian  Einstein-Weyl  $3$-manifolds $(M, [g], \triangledown)$; and 
\item orientation-reversing diffeomorphisms $\psi : \CP_1\to \CP_1$.  
\end{itemize}
\end{main} 
The prototypical example of such an 
  Einstein-Weyl  manifold is 
 three-dimensional de Sitter space 
 $SL(2,\CC)/SL(2, \RR)$; and the 
 diffeomorphism  $\CP_1 \to  \CP_1$ corresponding to this prototype is 
 the 
antipodal map 
\begin{equation}
\label{auntie}
[z_0 : z_1]\longmapsto [-\overline{z}_1 :  \overline{z}_0 ].
\end{equation}
In the present  context, 
two Einstein-Weyl structures are to be considered  identical 
if they are related by a connection-preserving conformal isometry.
On the other hand, 
  two orientation-reversing diffeomorphisms $\psi_1, \psi_2: \CP_1\to \CP_1$
are considered the same   iff 
$$\psi_1 = \varphi \circ
\psi_2 \circ 
\phi^{-1}$$
for   M\"obius transformations 
$\varphi , \phi  \in PSL(2, \CC )$ of the domain and range. 

In one direction, a direct geometrical interpretation of this  correspondence can
be described in terms of
{\em scattering maps}. We will thus  begin by explaining how to associate such a
map $\psi : \CP_1\to \CP_1$ to any Einstein-Weyl $3$-manifold $(M, [g], \triangledown)$  
satisfying the stated hypotheses.

 Recall that a conformal Lorentzian $n$-manifold $(M, [g])$ 
is said to be {\em space-time
oriented} if the structure group of its tangent bundle has been 
reduced to the identity  component $SO^\uparrow (1,n-1)\times \RR^+$ of 
the conformal Lorentz group;  in particular, this equips $M$ with a time orientation, allowing
one to determine whether a  time-like vector  is
past-  or future-pointing. A time-oriented Lorentz manifold
 $(M, [g])$  is  called {\em globally hyperbolic} if 
it contains a Cauchy surface, meaning a space-like hypersurface $\Sigma$ which meets
every endless  time-like or null curve exactly once; when this happens, there
is then  a diffeomorphism $M \approx \Sigma \times \RR$ such that 
every level set $\Sigma \times \{ t\}$ is a Cauchy surface    \cite{geroch,pentech}. 
We generalize the usual notion of conformal compactness for  Einstein manifolds 
\cite{grahlee,penbat} as in \cite{lebsdhg}, 
declaring that an Einstein-Weyl space $(M, [g], \triangledown)$ is 
{\em conformally compact} iff 
\begin{itemize}
\item there is a smooth compact pseudo-Riemannian
manifold-with-boundary $(X,\hat{g})$;
\item there is a diffeomorphism $\Phi: M \stackrel{\approx}{\longrightarrow} X - \partial X$
with   $[\Phi^*\hat{g}]= [g]$; 
\item the induced boundary metric  $\hat{g}|_{\partial X}$ is everywhere non-degenerate; 
and 
\item if $u$ is a non-degenerate defining function for $\partial X$, and 
if $\alpha$ is the  singular $1$-form   defined by 
$\triangledown \hat{g} = \alpha \otimes  \hat{g}$, then 
$\alpha - 2u^{-1} du$
is a smooth  $1$-form  
on $X$ with $d\alpha |_{\partial X}=0$.
\end{itemize}
If $(M,[g])$ is globally hyperbolic, then $\partial X \subset X$ is automatically  space-like,
and $X\approx \Sigma \times [0,1]$ for some compact $(n-1)$-manifold
$\Sigma$. If  $n=3$ and  $M$ is space-time oriented, 
it then follows  that  $M$  
 is  diffeomorphic to $\Sigma \times \RR$ for some compact oriented surface $\Sigma$, and   
 can  be conformally 
compactified by adding two copies ${\mathscr I}_-$ and  ${\mathscr I}_+$
of $\Sigma$, henceforth called {\em past} and 
{\em future  infinity}.

\begin{lem} \label{one}
 For an Einstein-Weyl manifold  as above, 
let $p\in {\mathscr I}_-$ be any point of past infinity. Then 
all the null geodesics emanating from  $p$  refocus 
at a unique point $q\in {\mathscr I}_+$. 
Moreover, ${\mathscr I}_-$ and  ${\mathscr I}_+$ are necessarily diffeomorphic to $S^2$. 
\end{lem}
\begin{proof}
First notice that every  null geodesic in $(M,[g])$
can be extended to a null geodesic in $(X, [\hat{g}])$, and therefore  has 
 unique past and  future end-points on ${\mathscr I}_-$ 
and ${\mathscr I}_+$, respectively. Indeed, the null geodesics meeting
${\mathscr I}_\pm$ can be parameterized by the  unit tangent bundle 
$UT{\mathscr I}_\pm$, and following these backwards or forwards 
to a Cauchy surface  $\Sigma\subset M$ then gives a local diffeomorphism  
$UT{\mathscr I}_\pm\to UT\Sigma$ which must be onto, since the image is both 
open and compact. 
 Using the conformal invariance of null geodesics, 
together 
with the fact that they are also geodesics of the Weyl connection $\triangledown$,
the hypothesized form of the singularity of $\triangledown$ at $\partial X$ 
 thus implies that  the affine parameter $t$ 
on any inextendible null $\triangledown$-geodesic must range over all of $\RR$, and 
must approach $- \infty$ (respectively,  $+\infty$)  exactly when the geodesic approaches 
its ideal  past (respectively,  future) end-point  
$p\in {\mathscr I}_-$  (respectively,  $q\in {\mathscr I}_+$).

This said,  let $\gamma : \RR \to M$ now be a $\triangledown$-affinely parameterized  null geodesic 
in $(M,[g])$. Let $\wee$
be a $\triangledown$-parallel non-zero space-like vector field along $\gamma$
which is orthogonal to the null vector field $\bfv:=\gamma^\prime$. 
Now observe that, because  
$\triangledown$ preserves the conformal structure, 
$${\mathcal R}\in \Lambda^2  \otimes (\Lambda^2 \oplus \RR g) $$
after index lowering. Hence 
 (\ref{hew})  
 implies that  ${\mathcal R}_{\bfv\wee}\bfv$ is $g$-orthogonal to  both 
$\bfv$ and $\wee$, and  is therefore a multiple of $\bfv$. 
It follows that there is a $\triangledown$-Jacobi field $\tilde{\wee}$
along $\gamma$ with $\tilde{\wee} \equiv \wee \bmod \bfv$. Thus $\wee \bmod \bfv$ actually represents an
infinitesimal  variation of $\gamma$ through (unparameterized)  null geodesics.
However, we have assumed that $\triangledown \hat{g} = \alpha \otimes  \hat{g}$
where 
$\alpha - 2u^{-1} du$ is a smooth up to ${\mathscr I}_\pm$,  so $t\sim \pm \log u$
as $\gamma$ approaches $\partial X$,  
and  $\hat{g} (\wee, \wee) \to 0$ as our affine parameter $t \to \pm \infty$. 
Hence $\wee$  actually represents an infinitesimal variation of $\gamma$ through  
null geodesics
emanating from   $p\in {\mathscr I}_-$ and terminating at 
 $q\in {\mathscr I}_+$. It follows that an infinitesimal change of  the initial direction of a null 
 geodesic emanating from 
  $p\in {\mathscr I}_-$  results in no change 
 at all  of the terminal end-point $q\in {\mathscr I}_+$. Since the set of future-pointing null directions
 at $p$  forms a connected family ($\approx S^1$), we thus conclude 
that all of the null geodesics emanating from $p\in {\mathscr I}_-$
 terminate at the same point  $q\in {\mathscr I}_+$.
 
 The family of null geodesics emanating from an arbitrary point $p\in {\mathscr I}_-$
 and refocusing at some  $q\in {\mathscr I}_+$  now gives us an immersion of 
 $\RR \times S^1$ in  $M$, since the corresponding Jacobi fields $\tilde{\wee}$ are everywhere 
 non-zero; moreover, its tangent space is   null  at every point. 
 The closure of this surface in $X$, obtained by adding $p$ and $q$, 
 is  therefore a topological immersion of $S^2$. 
 By  truncating  and
 capping  off with 
   space-like disks slightly to the future of $p$ and slightly to the past 
   of $q\in X$, we may then 
  smooth this 
 object  out  
 into a differentiably immersed $2$-sphere ${\mathcal S}\looparrowright M$
 such  that $T{\mathcal S}$ is then everwhere space-like or null. 
   Following a time-like vector field then gives us an immersion 
   $\varpi: {\mathcal S} \to \Sigma$ onto a Cauchy surface. But since $T{\mathcal S}$
   is everywhere  non-timelike, 
   $\varpi$ is necessarily  a covering map;  and  since $\mathcal S$ is simply connected, 
   this  must actually be the universal cover of 
   $\Sigma$. But since $\Sigma$ is an oriented 
   surface, the deck transformations of $\varpi$ must  be orientation-preserving. 
    Since any orientation-preserving diffeomorphism of $S^2$ has
   a fixed point, it therefore follows that there cannot be any non-trivial deck transformations of 
   $\varpi$, and 
   $\varpi$ must therefore 
   be  a diffeomorphism. Hence  $\Sigma$ and 
   $ {\mathscr I}_\pm\approx \Sigma$ are  $2$-spheres,  as claimed. 
  \end{proof}
 
 Since the $2$-spheres ${\mathscr I}_\pm$  
inherit an  induced conformal
structure from $X$,  and since both are  oriented by the space-orientation of
$M$, each is  isomorphic to  
  $\CP_1$ in a manner unique up to a M\"obius transformation.
We may thus define the  {\em scattering map}
$\psi : \CP_1\to \CP_1$ associated with $(M, [g], \triangledown)$     to be the
function  corresponding to 
$${\mathscr I}_-\ni p \longmapsto q\in {\mathscr I}_+$$ via some arbitrary 
choice  of  oriented 
 conformal isomorphisms 
${\mathscr I}_\pm\cong \CP_1$. Since the scattering $p \mapsto q$  
can  locally be realized  by 
following a congruence of null geodesics from ${\mathscr I}_-$ to  
${\mathscr I}_+$, $\psi$  is automatically smooth; and since 
we may construct a smooth inverse for it by instead following 
null geodesics backwards from ${\mathscr I}_+$ to  
${\mathscr I}_-$, we see that $\psi$ is in fact always a diffeomorphism. 

The fact that $\psi$ is necessarily 
orientation-reversing results from $q\in {\mathscr I}_+$ being the first point 
null-conjugate to $p\in {\mathscr I}_-$ along
any  null geodesic $\gamma$ joining them. 
Indeed, letting $\bfv$ and $\wee$ be as in the proof of 
Lemma \ref{one} above, and letting $g$ be 
a choice of $g\in [g]$ which is parallel along 
$\gamma$, there are $\triangledown$-Jacobi fields 
$\wee_1$ and $\wee_2$ along $\gamma$ such that 
$\wee_1 \equiv t \wee \bmod \bfv$ and $g(\wee_2, \bfv ) = 1$. 
The corresponding sections $\tilde{\wee}_1$ and $\tilde{\wee}_2$  of
the normal bundle of $\gamma$ then  join infinitesimally separated null geodesics,
and so, by conformal invariance, extend to $X$.
The determinant of $T_p{\mathscr I}_-\to T_q{\mathscr I}_+$ is now
explicitly given by 
$$(\tilde{\wee}_1 \wedge \tilde{\wee}_2)|_p \longmapsto   (\tilde{\wee}_1 \wedge \tilde{\wee}_2)|_q, $$
and  since $\tilde{\wee}_1\wedge \tilde{\wee_2}$ 
changes sign exactly once, at $t=0$, it thus follows that 
 the Jacobian determinant of the scattering map is negative at $p$.
The constructed map $\psi$ is therefore an orientation-reversing 
diffeomorphism, as promised. 
\bigskip 

We now turn to the problem of 
 inverting the  correspondence; and it is precisely here that twistor ideas
will  finally come to the fore. 
 The graph of any orientation-reversing diffeomorphism $\psi : \CP_1\to \CP_1$
 is a {\em totally real }  $2$-sphere  $P \subset \CP_1\times \CP_1$. 
Our program is then  to construct the corresponding 
$3$-manifold $M= M_\psi$ as a family of holomorphic disks
in $Z=  \CP_1\times \CP_1$ with boundaries on $P\subset Z$. When
$\psi$ is the antipodal map (\ref{auntie}),  these disks 
 are explicitly given by 
$$\zeta \longmapsto ([ a\zeta +b : c \zeta +d]  , [-\overline{d}\zeta-\overline{c} :
 \overline{b}\zeta + \overline{a}]) $$
 as $\zeta$ ranges over the unit disk $|\zeta|\leq 1$ in $\CC$. Here, of course, 
each 
$$ \left[\begin{array}{cc}a & b \\c & d\end{array}\right]
\in SL(2, \CC )$$
represents a particular choice of parameterized disk. 
 Notice that the 
 the boundaries of these disks
are exactly the standard  round circles in $P \cong S^2$,
and that the moduli space $M$ of disks modulo reparameterizations  
is exactly de Sitter space $SL(2,\CC)/SL(2, \RR )$.

If we form the double of 
any of the above disks, its abstract  normal bundle, in the sense of \cite{lebrsb}, will be 
${\mathcal O}(2)$. Since $h^1 (\CP_1 , {\mathcal O}(2))=0$ and 
$h^0 (\CP_1 , {\mathcal O}(2))=3$,  it follows that that they are all Fredholm regular, 
that the $3$-parameter family which they form is locally complete,
and that the  family is stable under totally real deformations of $P$; 
cf.  McDuff-Salamon \cite[Appendix C]{mcsalt} for
a different approach to  this  deformation problem.  
In particular,  if the totally real 
$2$-sphere $P\subset Z$ is perturbed slightly, the family 
of disks holomorphic disks $(D^2, \partial D^2) \hookrightarrow
(Z,P)$ will survive. In order
to study large perturbations of $P$, however, 
 one needs something like an energy estimate. Fortunately, the needed bound is freely
available in the present context, as a consequence of  the following observation: 

\begin{lem} \label{two} 
Let $\psi : \CP_1\to \CP_1$ be any orientation-reversing diffeomorphism, 
and let  $P \subset \CP_1\times \CP_1$ be the graph of $\psi$. Then there is a
K\"ahler metric $h$ on   $Z=\CP_1\times \CP_1$  such that $P$ is Lagrangian 
with respect to the corresponding K\"ahler form $\omega$, and such that 
$\omega$ has  de Rham cohomology class  $[\omega ]= 2\pi c_1 (Z)\in H^2 (Z, \RR )$.
\end{lem}
\begin{proof}
Let $h_2$ be the standard metric on the unit 2-sphere, which we identify in the 
usual way with $\CP_1$. Then $h_2$ is K\"ahler, with K\"ahler form $\omega_2$. 
Set $\omega_1 = -\psi^*\omega_2$. Since  $\psi$ is an orientation-reversing
diffeomorphism, $\omega_1$ is a positive $2$-form on $\CP_1$, and 
is therefore the area form of a unique conformal rescaling $h_1$ of $h_2$. Moreover, 
the total area of $\CP_1$ with respect to either $h_1$ or $h_2$ is $4\pi$,
so that $[\omega_1]$ and $[\omega_2]$ both represent $2\pi c_1 (\CP_1)$.

Now consider the product metric
 $$h=h_1\oplus h_2 :=  \varpi_1^*h_1 + \varpi_2^*h_2$$
where $\varpi_j  : \CP_1 \times \CP_1 \to \CP_1$, $j=1,2$,  are the factor projections. 
Then $h$ is a K\"ahler metric with K\"ahler form $\omega = \omega_1 \oplus \omega_2 :=
 \varpi_1^*\omega_1 + \varpi_2^*\omega_2$. The restriction of 
 $\omega$ to the graph of $P$ is then given by 
 $$(\mbox{id} \times \psi )^* (\varpi_1^*\omega_1 + \varpi_2^*\omega_2)=
 \omega_1 + \psi^*\omega_2 = \omega_1 - \omega_1= 0$$
 so  $P$ is therefore Lagrangian. Finally, 
 $$[\omega]=  \varpi_1^*[\omega_1] + \varpi_2^*[\omega_2]
 =  \varpi_1^*2\pi c_1 (\CP_1)  + \varpi_2^*2\pi c_1 (\CP_1) =2\pi c_1 (\CP_1 \times \CP_1),$$
and the K\"ahler class $[\omega]$ is therefore independent of $\psi$, as claimed. 
\end{proof}

If $\psi : \CP_1 \to \CP_1$ is an orientation-reversing diffeomorphism, 
it has degree $-1$, and  its graph $P$ therefore has homology 
class $(1,-1)$ in $H_2 (\CP_1\times \CP_1 , \ZZ ) \cong \ZZ \oplus \ZZ$. 
The long exact sequence
$$\cdots \to H_2 (P)  \to H_2 (Z)  \to H_2 (Z,P)   \to H_1 (P) \to \cdots$$
therefore tells us that $H_2 (\CP_1\times \CP_1,P;\ZZ) \cong \ZZ$,
where the generator ${\mathbf a}$ can be realized as the image of either 
$\CP_1\times \{ pt\}$ or $\{ pt\}\times \CP_1$ in  $\CP_1\times \CP_1$. 
As indicated by  the de Sitter space example discussed above, the 
holomorphic disks $(D^2, \partial D^2) \to
(Z,P)$ relevant to our problem  are precisely those which belong to
this generating class ${\mathbf a}\in H^2 (Z,P)$. With respect to the K\"ahler
metric $h$ of Lemma \ref{two}, these all have area $4\pi$. 

\begin{lem} \label{three} 
Let $\psi : \CP_1 \to \CP_1$ be any orientation reversing diffeomorphism, 
let $P\subset \CP_1\times \CP_1$ be its graph, 
and let  $F: (D, \partial D) \to
(Z,P)$ be any holomorphic disk representing the generator 
 $${\mathbf a}\in H_2 (\CP_1\times \CP_1,P;\ZZ) \cong \ZZ.$$ Then 
 $F$ is a holomorphic embedding,   is smooth up to the boundary, 
and  sends the interior of
 $D$ to the complement of $P$. Moreover,  
 $F (\partial D )$ is a smooth Jordan curve in $P \approx S^2$. 
 
 \medskip 

\noindent 
 Moreover, if $C$ is any holomorphic curve with boundary representing ${\mathbf a}$,
 then $C$ is non-singular, and  is either the image of a holomorphic disk
 as described above, or is a factor $\CP_1$ of $\CP_1 \times \CP_1$. 
\end{lem}
\begin{proof}
Consider the abstract oriented  $2$-sphere obtained by taking the double
$D \cup \overline{D}$, where the two copies of the disk are identified along the
boundary, $D$ is given the usual orientation coming from the
unit disk in $\CC$, and $\overline{D}$ is given the opposite orientation. 
Given $F$ as above, construct a  continuous map 
$\hat{F} : D \cup \overline{D}\to \CP_1$ by
$$
\hat{F} (z) = \left\{ 
\begin{array}{ll}
     \varpi_1 \circ F (z)   & \mbox{ if }  z\in D,   \\
      \psi^{-1}\circ   \varpi_2 \circ F (z)   &    \mbox{ if }  z\in \overline{D}.
\end{array}
\right.
$$
Since we have assumed that that $\psi$ is smooth, so is $P$, and 
 the holomorphic map $F$ is therefore  smooth up to the boundary  \cite{chirka}. 
 Hence $\hat{F}$ is 
actually smooth when  restricted to either $D$ or $\overline{D}$,
and is moreover orientation-preserving at every regular point.

Now, equipping $\CP_1$ and $\CP_1\times \CP_1$ respectively with the K\"ahler forms
 $\omega_1$ and $\omega$ of the proof of 
 Lemma \ref{two}, we have 
$$\int_{D \cup \overline{D}}\hat{F}^*\omega_1=
\int_{D}{F}^*(\varpi_1^*\omega_1+ \varpi_2^* \omega_2) = \int_{D}{F}^*\omega
= 4\pi = \int_{\CP_1}\omega_1$$
and it therefore follows that $\hat{F}$ has degree $1$. Moreover, since 
$\hat{F}$ is an orientation-preserving map at each of its regular points, 
the inverse image of any regular value must consist of exactly one point.
We therefore conclude that $F(D)$
is the graph of a Riemann-mapping biholomorphism 
 between two simply-connected domains with smooth boundary in 
 $\CP_1$. 
 
 If $C$ is a singular Riemann surface in $Z$ with boundary on 
 $P$, and if $[C ] = {\mathbf a}$ in relative homology, 
 we may first break $C$ up into its irreducible components, 
 and then observe that, because  each component must have area $\leq 4\pi$,  
  $C$ can in fact have only one irreducible component, with 
  multiplicity $1$, and so must be irreducible. If $C$ has non-empty boundary, 
 a straightforward generalization of the above doubling argument 
 shows  that $C \cup \overline{C}$ must be 
 homeomorphic to $\CP_1$, and it therefore follows that $C$ is
 one of the disks we have already analyzed. 
 On the other hand, if $C$ has no boundary, the fact that it must 
 once again have area $4\pi$ immediately implies that it must be a factor
 $\CP_1$. 
 \end{proof}

If $P\subset \CP_1\times \CP_1$ is the  graph of 
any orientation reversing diffeomorphism
$\CP_1 \to \CP_1$, 
Gromov's compactness theorem \cite[Theorem 4.6.1]{mcsalt} therefore tells us  that
the space of all embedded holomorphic disks $(D, \partial D) \subset 
(Z,P)$ with  $[(D, \partial D)]= {\mathbf a}\in H_2 (Z,P)$, together with the collection of
lines of the form $\CP_1\times \{ pt\}$ and $\{ pt\}\times \CP_1$,
forms a compact topological space. On the other hand, Lemma \ref{three}
also tells is than any such disk is 
diffeomorphically conjugate to one of the disks in our de Sitter example, 
and therefore has the same Maslov index. Since such a disk
is therefore \cite{lebrsb,mcsalt} Fredholm regular, deformation theory 
implies  that the space
of such disks is a smooth $3$-manifold --- although, {\em a priori},
it might still be  either empty or  disconnected. 

To show that the moduli space $M_\psi$ of disks is non-empty for 
each $\psi$, we first recall that the space of  orientation-reversing
diffeomorphisms of $S^2$ is connected, and that our deformation theory
tells us that  the set of $\psi$
for which disks {\em do} exist is open.
 It is therefore  tempting to 
just try to invoke Gromov compactness to show that it is also closed. 
However,  we know that sequences of disks can degenerate into factor $\CP_1$'s, 
and we must  control this degeneration in order to guarantee that 
we actually get  a disk in the limit. We do so by introducing the function 
$\Omega : M_\psi\to (0,4\pi) $ on the moduli space of disks which assigns 
to any disk $F: (D^2, \partial D^2) \to
(Z,P)$ the number $\Omega (F) = \int_{D^2} F^*\varpi_1^* \omega_2$, where $\omega_2$ is the
area form on the second factor. Since a sequence of disks $F_j$ can  degenerate
to $\CP_1\times \{ pt\}$ only if $\Omega (F_j) \to 0$, and can  degenerate to
 $\{ pt\}\times \CP_1$ only if $\Omega (F_j) \to 4\pi$, it follows that $\Omega$ is a proper map on
 the moduli space of disks. Taking limits as we vary $\psi$
 thus shows that, for every $\psi$,   the subset  $\Omega^{-1} (t) \subset M_\psi$ is 
  non-empty and compact 
 for each  $t\in (0, 4\pi )$; moreover, by Sard's theorem, 
 this level set  in any given $M_\psi$ is a smooth  compact surface for almost every $t\in (0, 4\pi )$.

 Given a choice of $\psi$, we  now turn to the construction of an Einstein-Weyl structure on 
 the moduli space $M=M_\psi$. Our first step in this direction is the construction
 of a conformal structure $[g]=[g]_\psi$ on $M$. To do this, we first
observe that our deformation theory   \cite{lebrsb} allows us to identify
 the tangent space of $M$ at a given disk
 $(D, \partial D) \subset (Z, P)$  with the
 space of  those holomorphic sections of the  normal bundle $(TZ|_D)/TD$ of $D\subset Z$
 for  which the boundary values are sections of 
 of the normal bundle $(TP|_{\partial D})/T\partial D$ of $\partial D\subset P$. 
 Because the Maslov index of the normal bundle of $D$ is $2$, this space
 is perfectly modeled by the space ${\mathfrak s \mathfrak l}(2, \RR)$
 of infinitesimal M\"obius transformations of the disk. 
 Consequently, up to rescalings,  each tangent space of $M$ carries a natural 
 Lorentz metric modeled on the 
 Killing form on  ${\mathfrak s \mathfrak l}(2, \RR)$. 
 The trichotomy of vectors into space-like, null, and time-like thus corresponds
 to the classification of   infinitesimal M\"obius transformations of the disk
 as hyperbolic, parabolic, or elliptic. A space-like vector on $M$
 therefore corresponds to an infinitesimal variation of the corresponding disk with two distinct zeroes
 on $\partial D$; a null vector corresponds to an infinitesimal variation with 
 a repeated zero on  $\partial D$; and a  time-like vector corresponds to 
 an infinitesimal variation with a single   zero in the interior of $D$,  but  none along $\partial D$.

In particular, a 
  one-parameter family of disks is time-like with respect to this conformal structure $[g]$ 
  iff the associated variation normal vector fields have no zeroes along the boundaries.
 By Lemma \ref{three}, this in particular means that the corresponding
 Jordan curves are nested, so as not to overlap, moving inward or outward at
 a non-zero rate everywhere. Projecting these disks
 into $\CP_1$ by $\varpi_2$, we thus see that the function 
 $\Omega: M\to (0, 4\pi )$,   which measures the area of these projections, 
 has non-vanishing derivative along any time-like curve.
 Consequently, $M$ carries a  time orientation for which $\Omega$ 
 is increasing along all future-pointing time-like curves. More strikingly, 
  the proper function 
 $\Omega
 : M\to (0,4\pi )$ therefore has no critical points,
 and integration of a  time-like vector field $v$ with $v\Omega =1$ thus gives us a
 diffeomorphism $M\approx \Sigma \times (0,4\pi )$, where $\Sigma$
 is the smooth compact 
 surface $\Omega^{-1}(0)$, and where  $\Omega$ has now become  factor projection to 
 the second factor 
 $(0,4\pi)$ of the product. 
 Since  the compactness of $\Sigma$ guarantees that any endless time-like or null curve
 in $M$ must  have  $\Omega \searrow 0$ at one extremity and $\Omega \nearrow 4\pi $
 at the other, the  intermediate value theorem forces 
 the level set $\Sigma \times \{ 0\}$ to meet every such curve.
Thus $(M,[g])$ has   a Cauchy surface, and   is therefore globally hyperbolic.

Now,  following any past-endless time-like curve backwards
in $M$  corresponds to a nested family of disks
in the second-factor $\CP_1$ with unique point  as intersection, since we know
that the Gromov limit  of the corresponding disks in $Z$ is one of the lines
$ \CP_1\times  \{ pt \} $. This gives us a unique limit point on $\CP_1$. 
Following the time-like factors of $M\approx \Sigma \times (0,4\pi )$ backwards
therefore gives us homeomorphism from $\Sigma$ to $\CP_1$,
and in particular giving $\Sigma \approx S^2$ a preferred orientation. 
Thus  $M\approx S^2 \times \RR$ is  
space-time oriented,
as promised.

We now endow $M$ with a Weyl connection $\triangledown$ compatible with $[g]$. 
To do so, we begin by describing its geodesics. The time-like geodesics
are the family of holomorphic disks through a given point
$x\in Z- P$.  The null geodesics are the families of  disks passing through a given 
point $x\in P$ with specified tangent. (These really are null geodesics of $[g]$,
because \cite{pentech}  they belong to the boundary of the future of a point, without
actually entering its future.) 
Finally, the space-like geodesics are
the families of disks passing though a pair of distinct points $x\neq y$ of $P$. 
Of course,  is not immediately obvious that there really {\em is} a connection
which has  these curves as its geodesics, but Fuminori Nakata has 
 written out a careful proof \cite[Proposition 8.1]{nakata} of this crucial fact
by  methods  based on  our   study of Zoll surfaces 
   \cite[Theorem 4.7]{lmzoll}. 
 Requiring the connection $\triangledown$  also be  $[g]$-compatible
 then specifies it completely. The fact that this connection is Einstein-Weyl then 
 follows \cite{hitproj}  from the fact that every null geodesic belongs to a totally 
 geodesic null hypersurface --- namely, the family of disks passing through 
 a given point of $P$.

However, we prefer a second  method of proving the existence of  
 $([g],\triangledown )$, as this   will turn out to be 
  better adapted to analyzing the geometry at infinity. 
 First notice that each of our disks $D$ has a canonical lift to the projectivized
 tangent bundle ${\mathbb P}T Z= {\mathbb P}T(\CP_1\times \CP_1)$,
given by its tangent bundle $TD\subset TZ$; moreover, 
 this lift is automatically Legendrian with respect to the tautological 
 contact structure on ${\mathbb P}T Z = {\mathbb P}T^* Z$. Moreover,
 since each  disk $D$ is the graph of a biholomorphism between 
 simply connected domains in $\CP_1$, its tangent line $TD$ is 
 never  vertical  or horizontal, and therefore may be viewed as a non-zero element 
 of the line bundle ${\mathcal O}(-2, 2)$ of linear maps
between the tangent space of the first and second factors. However, 
the complement of the zero section in ${\mathcal O}(-1, 1)$ 
is the universal cover of  the complement of the zero section in 
 ${\mathcal O}(-2, 2)$. Moreover, the complement of the  zero section
 in ${\mathcal O}(-1, 1)$ is exact $\CP_3$ minus two skew lines $L_1$ and $L_2$,
the projection to $\CP_1 \times\CP_1$ being given by 
 $[z^1:z^2:z^3:z^4] \mapsto ([z^1:z^2] , [z^3:z^4])$ if we choose our
 coordinates so that $L_1$ and $L_2$ are  given by 
 $z^1=z^2=0$ and $z^3=z^4=0$, respectively. We thus obtain
 two preferred liftings of each holomorphic disk $D$  to a 
 disk in $\CP_3$. Now  the inverse image of $P\subset \CP_1 \times \CP_1$ 
 is  a double cover of the projectivized tangent bundle $\RP TP$,  and so is 
 diffeomorphic to the unit tangent bundle $UTP$ of $P$ in a
 manner  which is completely canonical up to an overall choice of sign. 
 Making such a choice of sign once and for all, we then get a
  unique lift $\tilde{D}$ of each disk $D$ by requiring that $\partial\tilde{D}\subset UTP$
  is  the natural lift of the oriented Jordan curve curve  $\partial D \subset P$. 
  Each such disk $\tilde{D}\subset \CP_3$ then has its boundary on the 
  totally real submanifold $UTP\approx \RP^3$, and is Legendrian with respect
  to the  complex contact form   $\theta = z^1dz^2 - z^2 dz^1 + z^3 dz^4 - z^4 dz^3$.

  Each of these disks $\tilde{D}$ represents the generator of $H_2 (\CP_3, UTP )\cong \ZZ$,
so  the contact
  line bundle becomes ${\mathcal O}(2)$  on the abstract  double $\tilde{D}\cup \overline{\tilde{D}}$
  of any such disk. Because $\tilde{D}$ is  Legendrian, the normal bundle
  of its double is therefore $J^1{\mathcal O}(2) \cong {\mathcal O}(1) \oplus 
  {\mathcal O}(1)$. Moving these disks by the $S^1$ action $[z^1: z^2 : z^3 : z^4]
  \mapsto [z^1: z^2 : e^{i\theta}z^3 : e^{i\theta}z^4]$, we thus obtain a 
  $4$-parameter family $M\times S^1$ of disks in $(\CP_3 , UTP)$ 
  with doubled normal bundle ${\mathcal O}(1) \oplus 
  {\mathcal O}(1)$. It therefore follows  \cite[Prop. 10.1]{lmfrei} that 
  $M\times S^1$ carries an $S^1$-invariant self-dual conformal structure
  of signature $(++--)$, and the quotient of this geometry  by 
  $S^1$ therefore \cite{caldped,jonestod} gives us  the promised  Lorentzian Einstein-Weyl structure
  on $M$. 
  
  Now orbits of the  the $S^1$-action on $UTP\subset \CP_3$ are  the fibers 
  of the projection $UTP\to P$, and each of these circles is thus  the intersection
  of a fiber of ${\mathcal O}(-1,1) \to Z$ with $UTP$. 
 Each of these circles is therefore  contained in a unique 
  projective line $L$ in  $\CP_3$ which meets both of the skew lines $L_1$ and $L_2$. 
Such a circle then divides $L$  into two disks,  one of which 
 meets only $L_1$,  while the other meets only $L_2$. 
We thus obtain 
 two $S^2$-families of essentially explicit holomorphic disks $(\Delta , \partial 
\Delta )\subset (\CP_3 , UTP )$. 
For any such disk $\Delta$, 
 the abstract double  $\Delta\cup \overline{\Delta}$ is concretely
realized as  $L \subset \CP_3$, and  the double  therefore has
abstract  normal bundle ${\mathcal O}(1) \oplus {\mathcal O}(1)$.
Deformation theory \cite{glob,lebrsb} therefore gives us a
$4$-parameter family of disks in $(\CP_3, UTP)$ 
near any such $\Delta$. However, the infinitesimal
variations of such a disk which continue to
meet the relevant line $L_1$ or $L_2$ form a linear subspace
of codimension $2$. It follows that, aside from our initial
$2$-spheres of essentially explicit disks, all the other disks $\Delta^\prime$ 
of this family are disjoint from $L_1$ and $L_2$, and so are
 contained in ${\mathcal O}(-1,1)$ minus the zero section. 
 Projecting any such disk $\Delta^\prime$ into  $\CP_1 \times \CP_1$
 then gives  us a
holomorphic disk with boundary on $P$. However, since 
 the natural homomorphisms $H_2 (\CP_3 -(L_1\cup L_2), UTP) \to H_2 (\CP_3, UTP)$
and $H_2 (\CP_3 -(L_1\cup L_2), UTP)\to H_2 (\CP_1\times\CP_1, P)$
are both   isomorphisms, the images of  any of  these disks represents
the generator of $H_2 (Z,P)$. By Lemma \ref{three}, the projection
of any such $\Delta^\prime$  is therefore a disk $D$ of our original family. 
However, thinking of $\Delta^\prime$ as a section of ${\mathcal O}(-1,1)|_D$, 
we now see that it must be obtained from the canonical lift $\tilde{D}$ by 
multiplication by a holomorphic function on $D$ with unit modulus 
on $\partial D$; and since any such function is constant, it follows
that each $\Delta^\prime$  actually belongs to our previous family
$M\times S^1$. It follows that $M\times S^1$ can be compactified
by adding the two $S^2$-families of disks $\Delta$, and the
result
is now an $(S^2\times S^2)$-family of disks in $(\CP_3, UTP)$
with doubled normal bundle ${\mathcal O}(1) \oplus {\mathcal O}(1)$. 
The moduli space $Y\approx S^2 \times S^2$ therefore \cite{lmfrei} 
carries a split-signature self-dual metric which is invariant
under an $S^1$ action which just rotates one  $S^2$ factor about
an axis, while acting trivially on the other. The quotient 
$X = Y/S^1$ therefore \cite{lebsdhg} carries an induced conformal structure, 
and its interior $M$ carries an Einstein-Weyl structure 
whose singularity at $\partial X$ is exactly as described in our
definition of conformal compactness. Thus any 
orientation-reversing diffeomorphism $\psi: \CP_1 \to \CP_1$
determines a globally hyperbolic, 
conformally compact,  Lorentzian  Einstein-Weyl manifold, exactly as claimed.

  Of course, it  remains to check that our two constructions really are
  inverses of one another. In passing from the twistor-disk  construction 
  to the associated  scattering map, this is relatively straightforward, since 
  in the twistor  picture a null geodesic
  is just the family of disks with boundaries tangent to 
  a fixed element of $UTP$. The opposite direction is rather more subtle, 
  but just comes down to the fact that the mini-twistor space of a
  Lorentzian Einstein-Weyl space can be constructed, in analogy
  with \cite{hitmini} and \cite{lmzoll}, by starting with the space of future-pointing time-like
  geodesics, equipped with complex structure obtained by $90^\circ$ rotation of
  Jacobi fields,  and then adjoining a quotient of the space
  of null geodesics that represents the space of totally 
  geodesic  null surfaces. Details are left to the interested reader.

 We conclude by mentioning several surprising features of the 
 present correspondence  that might offer 
 fruitful directions for further investigation. First, 
the availability of  a K\"ahler metric for which $P$ is Lagrangian 
allowed us to prove  results here which are markedly sharper than 
those currently known  for  our earlier twistor-disk  constructions \cite{lmzoll,lmfrei}, and 
we therefore  wonder if some of our previous results could be sharpened 
by means of related techniques. Second, we find it interesting that,
{\em a posteriori},  any divergent 
sequence of   lifted disks $\tilde{D}_j\subset {\mathbb P}TZ$ must 
 subconverge to a singular curve consisting of  the Legendrian lift of a 
  $\CP_1$ and   half of a fiber of 
${\mathbb P}TZ \to Z$; and   we wonder if this might instead be shown more directly by 
some Gromov-type compactness argument. Third, the  Legendrian nature of 
the lifted disks $\tilde{D}$ means that the associated  embedding
 $M\hookrightarrow Y$ is actually umbilic, and that \cite{lebhspace} the complement of
a suitable hypersurface in $Y$ must therefore carry a  self-dual Einstein metric 
that one would like  to understand  more explicitly. 
Finally, as was   brought to our attention by
   Nakata's independent
investigation \cite{nakata},  
 the Einstein-Weyl 
structures constructed  here are all  actually 
{\em space-like Zoll}, in the sense that  their
space-like geodesics are all simple closed curves, and 
we are  intrigued by 
Nakata's very natural problem of classifying  all  possibles Einstein-Weyl $3$-manifolds with 
this 
property --- are they all given  by  the present construction? 
But  we  will  just leave these unresolved issues  for  the interested reader to ponder,
and simply   hope that  some of them may 
 eventually lead to interesting new avenues of  research.

  \end{document}